\documentclass[aps,prl,preprint,superscriptaddress]{revtex4-1}
\usepackage{graphicx}
\usepackage{dcolumn}
\usepackage{bm}
\usepackage{color}
\usepackage{amsmath}
\usepackage{natbib}
\begin{document}

%\preprint{}

\title{On the bursting of gene products}


\author{Romain Yvinec}
\email{yvinec@math.univ-lyon1.fr}
\affiliation{Institut Camille Jordan UMR 5208
Universit\'{e} Claude Bernard Lyon 1
43 Bd du 11 novembre 1918
69622 Villeurbanne cedex
France
}
\author{Alexandre F. Ramos}
\email{alex.ramos@usp.br}
\affiliation{Escola de Artes, Ci\^encias e Humanidades, Universidade
S\~ao Paulo, Av. Arlindo B\'ettio, 1000, CEP. 03828-000,
S\~ao Paulo, SP, Brazil
}
\thanks{This research was supported by the Natural Sciences and Research Council of Canada and the mobility fellowship CMIRA EXPLORA’ DOC 2010 of the r\'egion Rh\^one-Alpes (France). The work was carried out while both authors were visiting the Centre for Applied Mathematics in Bioscience and Medicine (CAMBAM). AFR and RY thanks CAMBAM members for warm hospitality and fruitful intelectual environment and specially M.C Mackey for carefully reading the manuscript. AFR thanks JEM Hornos for introduction into the field and partial support from Programa de Incentivo a Jovens Docentes of USP. RY thanks M. Santill\'an Zer\'on for useful comments.}
%\affiliation{}


\begin{abstract}
In this article we demonstrate that the so-called bursting production of molecular species during
gene expression may be an artifact caused by low time resolution in experimental data collection 
and not an actual burst in production.
We reach this conclusion through an analysis of a
two-stage and binary model for gene expression, and
demonstrate that in the limit when mRNA degradation is
much faster than protein degradation they are equivalent.
The negative binomial distribution
is shown to be a limiting case of the binary model for
fast ``on to off'' state transitions and
high values of the ratio between protein
synthesis and degradation rates.
The gene products population increases by unity but
multiple times in a time interval orders of magnitude
smaller than protein half-life or the precision
of the experimental apparatus employed in its detection.
This rare-and-fast one-by-one protein synthesis has been
interpreted as bursting.
\end{abstract}


\maketitle %must follow title, authors, abstract, \pacs, and \keywords

Understanding the origin of fluctuations at the single cell level and
how organisms deal them to guarantee both developmental viability and
evolutionary adaptation to a constantly changing environment
conditions is a challenge of the post-genomic era \cite{Kaern2005,
  Raser2005}.  Often, stochasticity at the single cell level is due to
the presence of biochemical reactants in low copy number inside the
cell \cite{delbruck40} and heterogeneous spatial distribution
\cite{Kuthan2001}.  Experimental techniques to investigate these
phenomena have been greatly enhanced by the use of fluorescent
molecules and technology to track the spatial and temporal behavior of
individual molecules \cite{Xie2008, Raj2009}.  Despite the striking
nature of the data these techniques provide, these advances do not
necessarily give the full picture of the dynamics of events at the
single cell level such as transcription and translation.  One example
is the measurement of the bursting production of molecules, defined as
an incremental increase in mRNA or protein number greater than one at
a given time.  Bursting is often held to be the {\it usual} mechanism
for the synthesis of gene products \cite{Xie2008, Raj2009}. As we show
here, the inference of bursting molecular production may be an
artifact due to a lack of sufficiently fine temporal resolution in
experimental data.  Also, from a modeling perspective the inference of
bursting may be flawed due to a reliance on the shape of stationary
probability distributions rather than an analysis of the underlying
dynamical processes giving rise to them.

The experimental observation of these jumps in molecular numbers has
motivated several models for the prediction and fitting of observed
data, {\em e.g.} by employing a Langevin approach (continuous) or the
master equation (discrete).  In the continuous case, stationary gamma
distributions for molecular concentrations are predicted along with
discontinuous trajectories for the corresponding stochastic process
\cite{Friedman2006, Mackey2011}.  In the discrete framework, bursts
appear in models for gene expression with two stochastic variables (so
called two-stage models with mRNA and protein), in the limit where the
mRNA degradation rate is much larger that the degradation rate for
protein (a common experimental finding).  In these cases, the model
predicted probability distributions are well described by a negative
binomial probability distribution \cite{Shahrezaei2008} and
simulations exhibit temporal bursting in protein numbers
\cite{Raj2009}.

In this article we use a discrete modeling framework to show that the
bursting limit actually corresponds to a particular regime of a model
of a switching gene between ``on'' and ``off'' states with a {\it one
  step variation} in the stochastic variable corresponding to protein
numbers \cite{Peccoud1995, sasai03, Hornos2005}.  The model we develop
here is an approximation to a model which is integrable in both the
stationary \cite{Hornos2005, innocentini07} and time-dependent regimes
\cite{Iyer-Biswas2009, ramos11} with symmetries underlying the
existence of analytical solution, and whose biological implications
have been explored elsewhere \cite{Ramos2007, ramos10}.

The steady state solution for the binary model, in the limit of rapid
transitions from the ``on'' to the ``off'' state, approaches a
negative binomial distribution that also describes the bursting model.
Simulation of the binary model shows apparent bursting, but examined
on a finer time scale reveals the protein numbers actually increase in
a unitary fashion.  This clearly suggests that the experimental
detection of bursting in gene product numbers may be due to lack of
temporal resolution in the data. Our results give clear guidelines for
the conditions on the transcription and translation processes for
artifactual bursts to appear.  As such, the modeling establishes
necessary conditions on the experimental temporal resolution necessary
to establish the existence of true bursting.

Gene expression is a cascade of first transcription to produce mRNA
followed by translation of that mRNA to protein.  Thus it makes sense
to take the system state variables to be the numbers of mRNA and
protein molecules in a given cell.  mRNA molecules are produced at a
rate that is dependent on the interaction between the RNA polymerase
and the promoter site. The number of mRNA molecules in the cell and
their interactions with ribosomes controls the protein synthesis.  In
this framework, {\it self regulation} is introduced in the model by
considering the transcription rate to be dependent on the number of
protein molecules produced.  It is worth noting that specific
regulation of the state of the promoter site (active or repressed) is
not taken into account on this model.

This picture for non-regulated genes has been treated
in the literature previously \cite{Thattai01}.

Let $m$ and $n$ denote, respectively, the number of mRNA and protein
molecules.  The probability of finding the system in a state $(m,n)$,
$m,n \geq 0$, at time $t$ is denoted by $P_{m,n}(t)$, while the
synthesis rates for mRNA and protein are denoted by $\mu^0_M$,
$\mu_M^1$ and $\nu_P$, and the corresponding degradation rates are
$\rho_M$ and $\rho_P$.  Then the evolution of the probability is
governed by a master equation for two coupled birth-death processes:
\begin{eqnarray}
\frac{dP_{m,n}}{dt} &=& (\mu^0_M+\mu^1_M n)(P_{m-1,n} - P_{m,n})\nonumber \\
                   &+& \nu_P m(P_{m,n-1}-P_{m,n}) \nonumber \\
                   &+& \rho_M[(m+1)P_{m+1, n} - P_{m,n}] \nonumber \\
                   &+& \rho_P[(n+1)P_{m,n+1} - nP_{m,n}], \label{eq1}
\end{eqnarray}
where we have assumed that the transcription rate is a function of the
protein number ($\mu_M^1 n$), indicating positive self regulation,
with the requirement that $\mu^0_M \neq 0.$ We have assumed a linear
dependence between the protein translation rate $(\nu_P)$ and the
number of available mRNA molecules in the cytoplasm.We have been
unable to construct an analytic solution to the complete system of
Eq. (\ref{eq1}).  However, exact quadrature is achieved in the
limiting case when the mRNA degradation rate is much greater than the
protein degradation rate ($\rho_M / \rho_P \gg 1$)
\cite{Shahrezaei2008} so the mRNA lifetime is quite short relative to
the protein lifetime.

That suggests scaling the model parameters by the protein lifetime
($\sim \rho_P^{-1} $), which results in the dimensionless quantities
\begin{equation} \label{dimlessparm}
\mu^0=\frac{\mu^0_M}{\rho_P}, \ \ \mu^1 = \frac{\mu^1_M}{\rho_P}, \ \ 
\gamma = \frac{\rho_M}{\rho_P}, \ \ \nu = \frac{\nu_P}{\rho_P}.
\end{equation}
The approximate Eq. (\ref{eq1}) becomes
\begin{eqnarray}
\frac{dP_{0,n}}{d\tau} &=& [(n+1)P_{0,n+1} - nP_{0,n}]  \nonumber \\
                      &-& (\mu^0 + \mu^1 n) P_{0,n}
                      + \gamma P_{1,n}, \label{eq2} \\
\frac{dP_{1,n}}{d\tau} &=& \nu (P_{1,n-1} - P_{1,n})
                        +  [(n+1)P_{1,n+1} - nP_{1,n}] \nonumber \\
                       &+& (\mu^0 + \mu^1 n) P_{0,n} - \gamma P_{1,n},
                          \label{eq3}
\end{eqnarray}
where we have introduced the dimensionless time $\tau = \rho_P t$
scale and the approximations $ P_{m,n} \sim 0, m\ge2$ and $(\mu^0 +
\mu^1 n)P_{1,n}/(\gamma P_{2,n}) \sim 1$, for all $n$.  Since our
simplifying assumption implies that the mRNA lifetime is short
relative to that of the protein, we would expect that mRNA
probabilities will be peaked around zero for $\mu^0, \mu^1$ of the
same order as $\rho_M$. This offers some justification for assuming
that Eqs. (\ref{eq2}) and (\ref{eq3}) are valid for describing gene
expression (Supplementary information).

Eqs. (\ref{eq2}) and (\ref{eq3}) have the same form as the master
equation for a binary gene with the state $(1,n)$ (or $(0,n)$) as the
active (or repressed) state of protein synthesis with rate $\nu$ (or
zero). The ``on-off'' switching rate is given in terms of the mRNA
degradation rate $\gamma$ and the ``off-on" transition depends on the
mRNA synthesis rates, $\mu^0$ (for external regulation) or $\mu^0 +
\mu^1 n$ (self regulation).  To write the solutions of the model
presented at the Eqs. (\ref{eq2}) and (\ref{eq3}), we define constants
$(a, b, \theta)$ as folows:
\begin{equation}
a = \frac{\mu^0}{1+\mu^1}, \ \ \ \
b = \frac{\mu^0+\gamma}{1+\mu^1}, \ \ \ \
\theta = \frac{1}{1+\mu^1},
\label{params}
\end{equation}
where external regulation is recovered by setting $\theta=1$.  For
simplicity we consider only the steady state solutions for
Eqs. (\ref{eq2}) and (\ref{eq3}), $P_{0,n}, P_{1,n}$, and the
probabilities for finding $n$ proteins inside the cell, $P_n= P_{0,n}
+ P_{1,n}$, namely:
\begin{eqnarray}
P_{0,n} &=&  \frac{b-a}{Cb}\frac{\nu^n}{n!} \frac{(a)_{n}}{(1+b)_{n}}
          {\rm M}(a+n, 1+b+n, -\nu \theta), \label{p0n} \\
P_{1,n} &=&  \frac{a}{Cb}\frac{\nu^n}{n!} \frac{(1+a)_{n}}{(1+b)_{n}}
          {\rm M}(1+a+n, 1+b+n, -\nu \theta), \label{p1n} \\
P_n    &=&  \frac{\nu^n}{Cn!} \frac{(a)_{n}}{(b)_{n}}
          {\rm M}(a+n, b+n, -\nu \theta), \label{phin}
\end{eqnarray}
where $M(a,b,z)$ denotes the Kummer M function \cite{abramowitz72} and
the normalization constant $$C={\rm M}(a, b, \nu(1-\theta)),$$ assures
conservation of probability
$\sum_{n=0}^{\infty}(P_{0,n}+P_{1,n})=1$. Note that for external
regulation, $C={\rm M}(a, b, 0) = 1$.

The denominator in Eq. (\ref{params}), $1+\mu^1,$ should be
interpreted as the total protein removal rate from the cytoplasm by
degradation plus transcription stimulus.  The constant $a$
characterizes the rate of spontaneous (basal) mRNA synthesis relative
to the protein removal rate and states a relation with the probability
for finding one mRNA ($p_1$), defined as $\sum_{n=0}^\infty P_{1,n}$,
namely
\begin{equation}
p_1 = \frac{a}{C b}\,
{\rm M}(a+1,b+1,\nu(1-\theta)),
\end{equation}
and, for the external regulating gene, $\theta=1$, it
implies $$p_1=\frac{a}{b}.$$ Phenomenologically, $b$ is a compound
relation between the rate for a cycle of mRNA synthesis-degradation
and the protein removal rate.  Its role in determining the statistics
of protein numbers is seen from the average number of protein
molecules, $\langle n \rangle = p_1 \nu $, and the variance relative
to the average (Fano factor), $\sigma^2/\langle n \rangle = (\langle
n^2 \rangle - \langle n \rangle ^ 2)/ \langle n \rangle$, given by
\begin{eqnarray}
\frac{\sigma^2}{\langle n \rangle} &=&
1 + \nu\frac{a+1}{b+1}
\frac{{\rm M}(a+2, b+2, \nu(1-\theta))}{{\rm M}(a+1, b+1, \nu(1-\theta))}
\nonumber \\
&-& \nu\frac{a}{b}
\frac{{\rm M}(a+1, b+1, \nu(1-\theta))}{{\rm M}(a, b, \nu(1-\theta))}. \label{fanobin}
\end{eqnarray}
For the case where $\theta=1$, we have
\begin{equation}\label{fanoext}
\frac{\sigma^2}{\langle n \rangle} = 1 + \frac{\nu}{b}\frac{1-a/b}{1+1/b}.
\end{equation}
In the limit $b\rightarrow +\infty$ (fast switching) and the
parameters $(a, \theta, \nu)$ are finite, Eq. (\ref{fanobin}) is equal
to one and the distribution of protein is Poissonian.

To get some intuition into this system, consider the steady state (or
equilibrium) of the Eqs. (\ref{eq2}) and (\ref{eq3}), when
probabilities are fixed with time, but the variables $(m,n)$ change in
time with the probabilities given by Eqs. (\ref{p0n}), (\ref{p1n}) and
(\ref{phin}).  $b$ gives the average time for the system to complete
one switching cycle -- {\em e.g.} from {\em off} to {\em on} and back
to {\em off} again.  Then the probabilities $p_1$ and $1-p_1$ are the
fractions of the total switching time that the system spends in the
active and inactive states respectively.

The transition from the dynamic to the stationary regime, as noted
previously \cite{innocentini07, ramos11}, has an approach to
equilibrium characterized by two of the time scales of the model,
$\rho_P^{-1}$ and $(\rho_P+\mu^1_M)^{-1}b^{-1}$.  The former is the
typical lifetime of the protein, whereas the second one is related to
the switching.  For rapid protein degradation, compared to the
switching, the steady state is achieved after the equilibrium of {\em
  on-off} transitions that are slow and can result in super Poisson
stationary distributions (eventually, bi-modality occurs with each
peak related to one state of the system).  On the other hand, when
protein degradation dominates, and there is fast switching, the
distributions are uni-modal. In that case, the gene regulatory
mechanism ({\em e.g.} if binary or constitutive) is indistinguishable
by simple protein counting.

This reasoning suggests that bursting would occur for systems with a
large value of $\nu$ and $p_1\sim0$.  Biologically, this would mean
that the mRNA number is mostly zero during an entire switching cycle.
For a $p_1$ fraction of that cycle, there will be one mRNA that is
rapidly translated (at rate $\nu$) and thus several unitary increments
in $n$ take place during a very short time.  This will appear to be a
single near-instantaneous increase in protein number by more than
one. A mechanism for a rapid increase in $n$ from one mRNA is the
binding of several ribosomes to the mRNA.

Mathematically, the negative binomial distribution is assumed to
describe a random variable characterizing a bursting process. We can
show (Supplementary information) that the negative binomial
distribution is a particular case of the probabilities of the
Eq. (\ref{phin}) at the limit of $b,\nu \rightarrow \infty$ with their
ratio
\begin{equation}\label{delta}
\delta=\frac{\nu}{b},
\end{equation} 
kept finite, namely:
\begin{equation}\label{negbin1}
P_n \rightarrow \frac{(a)_n}{n!}
\left( \frac{\delta}{1+\delta \theta} \right)^n
\left( \frac{1+\delta(\theta-1)}{1+\delta \theta} \right)^a,
\end{equation}
where $(a)_n=a (a+1) \dots (a+n-1), (a)_0=1$.  For the self-regulating
case, an approximate negative binomial distribution occurs for $\theta
\sim 1$, which implies weak induction of mRNA synthesis by proteins,
($\mu^1 << 1$).  For the externally regulated gene, $\theta=1$, and
the probabilities at equation above become:
\begin{equation}\label{negbin}
P_n \rightarrow \frac{(a)_n}{n!}
\left( \frac{\delta}{1+\delta } \right)^n
\left( \frac{1}{1+\delta \theta} \right)^a,
\end{equation}
that is the negative binomial distribution \cite{Shahrezaei2008}.

\begin{figure}
 \includegraphics[width=0.7\linewidth]{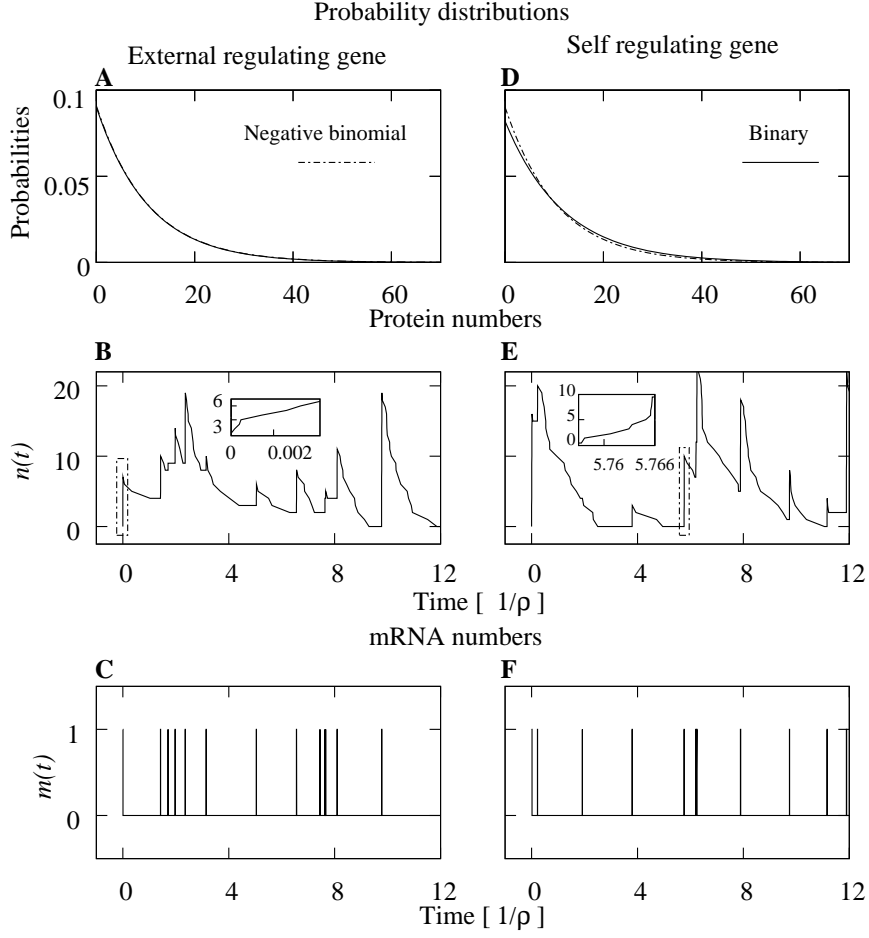}%
 \caption{ {\bf Steady state probability distributions and
     trajectories.}  The parameter values for the externally regulated
   gene are $\mu^0= 1, \gamma = 99$ and for the self regulating gene
   $\mu^0 = 100/99, \mu^1 = 1/99, \gamma = 100$.  The parameters $(a,
   b, \theta, \nu)$, for the binary models, are then $(1, 100, 1,
   1000)$ for the external and $(1, 100, 0.99, 1000)$ for the self
   regulating probabilities.  For the negative binomial distribution
   in FIGS. A and B we took $a=1$ and $\delta=\nu/b = 10$.  In
   FIG. \ref{figures}. B, E, C, and F, the trajectories for the binary
   model are shown in protein half-life time scale.  An expansion of
   the time scale in the region bounded by the dashed lines in
   FIG. \ref{figures}. B and E is shown in the insets where the time
   scales are magnified by $10^3$, the same as the ratio between
   protein synthesis and degradation rates. For the external and self
   regulating genes, the probabilities for finding one mRNA are,
   respectively, 0.01 and 0.011, the mean protein number $\langle n
   \rangle$ are, 10 and 11.08, and finally the Fano factor's values
   are 10.8 and 11.08.  The Fano factor value of the negative binomial
   distribution is 11, considering the same set of parameters, and is
   closely related to the one for the binary probabilities.}
  \label{figures}
 \end{figure}

We illustrate our results in FIG. \ref{figures} where the left hand
column is for the external regulated gene and the right hand column is
for the self regulating case. FIG. \ref{figures}.A and \ref{figures}.D
are the steady state probability distributions as obtained from the
expression of the binary (Eq. \ref{phin}) and negative binomial
(Eq. \ref{negbin}) models.  For the parameter values we choose,
inspection shows high agreement with the externally regulated while a
slight difference appears for the self-regulating gene.  The
corresponding trajectories are obtained from the binary model, with
the protein numbers shown in FIG. \ref{figures}.B and
\ref{figures}.E. The {\it apparent} bursts appear explicitly at
protein half-life time scale. However, an expansion of the time scale
reveals that the protein synthesis is occurring one by one.  Finally,
the corresponding mRNA number dynamics is shown in
FIG. \ref{figures}.C and \ref{figures}.F. As expected, it is switching
between very short time intervals with one mRNA and long intervals
with no mRNA.

Experimentally, the temporal resolution necessary for avoiding
anomalous bursting detection should be of the order of the average
time for translation of one protein, $\sim 1/\nu_P$. In what follows,
we shall consider the system approached in Ref. \cite{cai06} to show
an example of a measurement of apparent bursting. In their work, the
authors monitored the expression of the $\beta$-gal protein under the
control of the {\em lac} promoter. They have detected the occurrence
of burstings in protein numbers and measured the average bursting size
to be of $\sim 8$ proteins synthesized per burst. Our aim is to
calculate the protein synthesis rate of the system reported in
Ref. \cite{cai06} using the average bursting size calculated by the
authors. We shall employ our approximation at Eqs. (\ref{eq2}) and
(\ref{eq3}) and estimate the necessary time resolution for avoiding
the measurement of apparent bursting.

We start by setting the average bursting size in terms of the rates of
the Eq. (\ref{eq1}). In literature \cite{Shahrezaei2008}, the average
bursting size is usually given by the parameter $\delta$, of the
negative binomial probability distribution at the Eqs. (\ref{negbin1})
and (\ref{negbin}). Therefore, the protein synthesis rate, $\nu_P$ at
the Eq. (\ref{eq1}), can be written as a function of $\delta$ as
follows:
\begin{equation}\label{psynthesis}
\nu_P = \delta \, \rho_P \frac{\mu_M^0+\rho_M}{\rho_P + \mu^1_M},
\end{equation}
that is deduced by combining the Eqs. (\ref{delta}), (\ref{params}),
and (\ref{dimlessparm}).

To proceed calculating $\nu_P$, we set $\mu^1_M = 0$ -- and assume
external regulation -- since the {\em lac} promoter interacts with the
Lac repressor protein that {\em is not} encoded in the {\it lac
  operon} genes. Then, the expression for calculating the protein
synthesis rate at the Eq. (\ref{delta}) is reduced to $\nu_p = \delta
(\mu^0_M+\rho_M).$ The values of the remining unknown constants, mRNA
synthesis and degradation, are determined in terms of experimental
measures.

The mRNA degradation rate of the $\beta$-gal mRNA -- $\rho_M$ -- is
taken to be $\sim 0.1$ min$^{-1}$ based on data reported in
Ref. \cite{dickson81}.

We use the data provided in Ref. \cite{cai06} to estimate the mRNA
synthesis rate ($\mu^0_M$) at $ \sim 10^{-3}$ min$^{-1}$.  This number
is achieved dividing the average frequence of bursting of 0.16 per
cell cycle by the average period of a cell cycle, 145 min (both data
from Ref. \cite{cai06}). The bursts of proteins occur whenever one
mRNA arrives at the cytoplasm. Therefore, we can convert the average
bursting frequence into the average mRNA synthesis rate.

Based on the values of $\delta, \mu_M^0, \mu_M^1,$ and $\rho_M$ above
we compute the protein synthesis rate as $\nu_P \sim 1$ min$^{-1}$.
Thus, the time scale for the synthesis of one protein is $\sim 1$ min
or, $\ln{2}/\nu_P \sim 0.5$ min in case of exponential growth of the
protein population. Such time scales are smaller than the time
resolution of protein detection of 4 min reported in
Ref. \cite{cai06}. Intuitively, one might expect an average of at
least 4 protein synthesis during the temporal range of the
experimental resolution. The detection of an increase greater than one
in protein population should be interpreted as a bursting. However,
under the conditions we are considering in this manuscript, an
experimental time resolution of $\sim 1$ min would re-establish a
one-by-one protein population increase.

We also stress another aspect of experimental measurements: the
sampling time.  Usually, this is the time interval between two
gatherings of cells from the cell culture. In the experiments
presented in Ref. \cite{cai06}, the sampling time is 20 s.  Note that
it appears to be enough to detect the one-by-one protein
increment. However, the time resolution for protein detection is
obtained indirectly, from fluorescence measurements, which results a
4 min precision. Therefore, this is not enough to detect individual
proteins.

It is worth to discuss the phenomenology of the burst like behavior
reported in Ref. \cite{cai06}. Theoretically, it relates to the value
of the ratio $\nu = \nu_P/\rho_P$ that is a large number, $>>
1$. Then, whenever an mRNA appears in the cytoplasm, a plethora of
proteins is fastly synthesized while their degradation is very slow.
In a plot of the protein number versus time, that condition appears as
a fast increment in protein population, followed by a plateau, if the
experiment stands shorter than the protein half-life time.

In that case, a description of protein numbers inside the cell in
terms of the negative binomial distribution is appropriate. However,
it must be emphasized that the fitting of the measured histogram by a
known probability distribution does not imply the occurrence of the
stochastic process from where the distribution is derived. The
predictive power of a model based on such approach might be lowered.
In that sense, our discussion is an increment in the capability of
interpreting experimental data. That is done establishing the time
resolution necessary for experimental demonstration of bursting.

We establish two possible sceneries for the occurrence of real burst of
proteins. It is assumed experimental time resolution of the order of
$\sim 1/\nu_P$ and the measurement of a greater than one instantaneous
increment in protein population. The first scenery requires the
existence of only one mRNA in the cytoplasm. A possible mechanism to
underlie this effect is: multiple polipeptide chains present in
cytoplasm that were translated individually and start their functional
activity simultaneously. Note that that implies a delay between the
translation process and the protein folding. In our second bursting
scenery, there is abundant fast degradating mRNA's in the cytoplasm
that can be translated synchronously in an one-by-one fashion. Hence,
multiple proteins could be synthesized simultaneously in a time scale
of the order of $\sim 1/\nu_P$.

Our results also suggest a picture of gene expression where the
bursting (or burst-like) dynamics corresponds to one among other
possible behaviors. Different regimes of gene expression are possible
depending on the specific relations among the effective rates of the
reactions participating of a gene network. For example, the model at
Eq. (\ref{eq1}) has a precise biological interpretation. Its
approximation in terms of the binary model, Eqs. (\ref{eq2}) and
(\ref{eq3}), shows the utility of the ``on'' and ``off'' model for the
analysis of gene products synthesis. In terms of probability
distributions one expect that, besides the negative binomial,
the gene products should also satisfy the probabilities
based on Eqs. (\ref{p0n}), (\ref{p1n}) and (\ref{phin}). 

These probabilities satisfy the Eq. (\ref{eq1}) when the probability
to find more than one mRNA in the cytoplasm is neglegible. In this
sense, one might provide further approximations to the solution of the
Eq. (\ref{eq1}) for neglegible probability of detecting three, four
mRNA's in the cytoplasm, respectively, in terms of ternary, quaternary
models. Different insights in the workings of gene expression could
also be provided by this kind of approach.

In this manuscript we are proposing a theoretical framework, based on
the binary model to gene expression, that generalizes the negative
binomial distribution for the description of the stochasticity in gene
products number. The burst like behavior occurs in a well defined
regime, when the ratio between synthesis and degradation rates is of
the order of 10$^3$ and the synthesis of gene products very rare. As
we predict from our model, measurements aiming to detect the
one-by-one increments in gene products number must have temporal
resolution of the order of the synthesis rate ($1/\nu_P$), {\em e.g.}
in conditions reported in Ref. \cite{cai06} that would imply a time
resolution of $\sim 10-60$ s.

\bibliography{main_R_biblio}

\end{document}